[1994/12/01]% LaTeX date must December 1994 or later
\documentclass[draft]{article}
\usepackage{amsmath}
\usepackage{amsfonts}
\usepackage{amssymb}
\usepackage{euscript}
\input epsf
\usepackage{amsthm}
\newtheorem{thm}{Theorem}[section]
\newtheorem{cor}[thm]{Corollary}

\newtheorem{lem}[thm]{Lemma}

\theoremstyle{definition}

\theoremstyle{remark}

\begin{document}
\begin{center}

{\large \bf THE INEQUALITIES FOR \\SOME TYPES OF $q$-INTEGRALS}

\bigskip

{
 % PLEASE INSERT THE NAMES AND ADDRESSES OF THE AUTHORS HERE

\centerline{\bf Predrag M. Rajkovi\'c} \centerline{Depart. of
Mathematics, Faculty of Mechanical Engineering} \centerline{ {\it
e-mail}: {\tt pecar@masfak.ni.ac.yu}}

\smallskip

\centerline{\bf Sladjana D. Marinkovi\'c} \centerline{ Dep. of
Mathematics, Faculty of Electronic Engineering} \centerline{{\it
e-mail}: {\tt sladjana@elfak.ni.ac.yu}}

\smallskip

 \centerline{\bf Miomir S. Stankovi\'c}
\centerline{Department of Mathematics, Faculty of Occupational
Safety}
 \centerline{ {\it e-mail}: {\tt miomir.stankovic@gmail.com}}

\smallskip

\smallskip
\smallskip

\centerline{\bf University of Ni\v s, Serbia and Montenegro} }

\end{center}

% PLEASE INSERT THE ABSTRACT HERE

{\bf Abstract.} We discuss the inequalities for $q$-integrals
because of the fact that the inequalities can be very useful in
the future mathematical research. Since $q$-integral of a function
over an interval $[a,b]$ is defined by the difference of two
infinite sums, there a lot of unexpected troubles in analyzing
analogs of well-known integral inequalities. In this paper, we
will signify to some directions to exceed this problem.

\medskip

%33D60 Basic hypergeometric integrals and functions defined by them
%26D15 Inequalities for sums, series and integrals

\smallskip


\section{Introduction}

The integral inequalities can be used for the study of qualitative
and quantitative properties of integrals. In order to generalize
and spread the existing inequalities, we specify two ways to
overcome the problems which ensue from the general definition of
$q$-integral. The first one is the restriction  of the
$q$-integral over $[a,b]$ to a finite sum (see \cite{Gauchman}).
The second one is indicated in \cite{RSM} and it means
introduction the definition of the $q$-integral of the Riemann
type. At the start sections, we give all definitions of the
$q$-integrals, their correlations and pro\-per\-ties. In the other
sections, we elaborate the $q$-analogues of the well--known
inequalities in the integral calculus, as Chebyshev, Gr\"uss,
Hermite-Hadamard for all the types of the $q$-integrals. At last,
we give a few new inequalities which are valid only for some types
of the $q$-integrals.

In the fundamental books about $q$-calculus (for example, see
\cite{Gasper} and \cite{Hahn}),
 the $q$-integral of the function $f$ over the
interval $[0,b]$ is defined by
\begin{equation}
I_q(f;0,b)=\int_0^b f(x)d_q x =
b(1-q)\sum_{n=0}^{\infty}f(bq^n)q^n\quad (0<q<1).
\label{(1.1)}
\end{equation}
If $f$ is integrable over $[0,b]$, then
$$
\lim_{q\nearrow 1} I_q(f;0,b) = \int_0^b f(x)\ d x = I(f;0,b).
$$
Generally accepted definition for $q$-integral over an interval
$[a,b]$ is
\begin{equation}
I_q(f;a,b)=\int_a^b f(x)d_q x= \int_0^b f(x)d_q x - \int_0^a f(x)
d_q x\quad (0<q<1). \label{(1.2)}
\end{equation}
The values of such defined $q$-integrals of the polynomials have
very similar form to those in the standard integral calculus. So,
for example, we  it is valid
$$
\int_{a}^{b} x^n d_q x=\frac{b^{n+1}-a^{n+1}}{[n+1]_q},
$$
where
$$
[n]_q=\frac{1-q^n}{1-q}=1+q+\cdots+q^{n-1},\qquad n\in\mathbb N.
$$

But, the problems come when  the integrand $f$ is defined in
$[a,b]$ and it is not defined in $[0,a]$. Obviously this
definition cannot be applied on evaluating of the integrals of the
form
$$
\int_2^3 \ln (x-1)d_q x.
$$

\section{The $q$-integrals and correlations }

Let $a, b$ and $q$ be some real numbers such that
$0<a<b$ and $q\in (0,1)$.

Beside the $q$-integrals defined by (\ref{(1.1)}) and
(\ref{(1.2)}) we will consider two other types of the
$q$-integrals.

In the paper \cite{Gauchman}, H. Gauchman has introduced {\it the
restricted $q$-integral}
\begin{equation}
G_q(f;a,b)=\int_{a}^b f(x)\ d_q^{G} x=
b(1-q)\sum_{k=0}^{n-1}f(bq^k)q^k \quad (a=bq^n).
\label{(2.3)}
\end{equation}
Let us notice that lower bound of integral is $a=bq^n$, i.e. it is
tied by chosen $b$, $q$ and positive integer $n$.

In the paper \cite{RSM}, we have introduced {\it Riemann-type
$q$-integral} by
\begin{equation}
R_q(f;a,b)=\int_a^b f(t)d_q^{R}t=
(b-a)(1-q)\sum_{k=0}^{\infty}f\bigl(a+(b-a)q^k\bigr)q^k.
\label{(2.4)}
\end{equation}
This definition includes only point inside the interval of the
integration.

The different types of the $q$-integral defined by (\ref{(1.1)})--(\ref{(2.4)})
can be denoted in the unique way by $J_q(\,\cdot\,;a_{(J)},b)$,
where $J$ can be $G$,\ $I$\ or $R$. Interval of the integration
$E_{(J)}=[a_{(J)},b]$  of  $q$-integral $J_q(\,\cdot\,;a_{(J)},b)$
depends on its type:

 $a_{(G)}=bq^n$, $n\in\mathbb N$,
for $G_q(\,\cdot\,;a,b)$;

 $a_{(I)}=0$, for $I_q(\,\cdot\,;0,b)$;

  $a_{(I)} \in [0,b]$,  for $I_q(\,\cdot\,;a,b)$;

 $a_{(R)} \in [0,b]$,  for $R_q(\,\cdot\,;a,b)$.

\smallskip

We can say that a real function $f$ is {\it $q$-integrable} on
$[0,b]$ or $[a,b]$ if the series in (1) and (2) converge. In the
similar way, we say that $f$ is {\it $qR$-integrable} on $[a,b]$
if the  series in (4) converges.

From now on, it will be assumed that the function $f$ is
$q$-integrable on $[0,b]$ \ ($qR$-integrable on $[a,b]$) whenever
$I_q(f;0,b)$ or $I_q(f;a,b)$ \ ($R_q(f;a,b)$) appears in the
formula.

In this research it is convenient to define the operators \ \
$$
\aligned
&\widehat {}\ :f\mapsto \widehat  f,\qquad \widehat f(x)=f\bigl(a+(b-a)x\bigr),\\
&\widetilde {}\ :f\mapsto \widetilde f,\qquad  \widetilde f(x)=bf(bx)-af(ax),\\
&\breve{}\ :f\mapsto\breve f,\qquad \breve f(x)=f(bx)-f(ax),
\endaligned
$$
such that associate the functions defined on $[0,1]$  to the
function defined on $[a,b]$. Notice that, for $x\in[0,1]$, it is
\begin{equation}
\widehat{(fg)}(x)=\widehat f(x)\ \widehat g(x),\quad
\widetilde{(fg)}(x)=\frac 1{b-a}\bigl(\widetilde f(x) \widetilde g(x)-ab\
\breve f(x) \breve g(x)\bigr). \label{(2.10)}
\end{equation}

The correlations between the $q$-integrals defined by (\ref{(1.1)})--(\ref{(2.4)})
are given in the following lemma.

\begin{lem}\label{2.1}
If the real function $f$ is $q$-integrable on $[0,b]$ or
$qR$-integrable on $[a,b]$\ $(0<a<b)$, then it holds
\begin{eqnarray}
I_q(f;0,b)&=&\lim_{n\to\infty}G_q(f;bq^n,b), \label{(2.5)}\\
I_q(f;a,b)&=&I_q(\widetilde f;0,1),\quad \textrm{where}\quad
\widetilde f(x)=bf(bx)-af(ax), \label{(2.6)}\\
R_q(f;a,b)&=&(b-a)I_q(\widehat f;0,1), \quad\textrm{where}\quad
\widehat f(x)=f(a+(b-a)x), \label{(2.7)}
\end{eqnarray}
\end{lem}

\noindent{\it Proof.}  Since $G_q(f;bq^n,b)$ $(n\in\mathbb N)$ is
the
 partial sum of the series $I_q(f;0,b)$, the relation (\ref{(2.5)})
is evident.

The equalities (\ref{(2.6)}) and (\ref{(2.7)}) are valid because
of
$$
I_q(f;a,b)
=(1-q)\sum_{k=0}^{\infty}\bigl(bf(bq^k)-af(aq^k)\bigr)q^k
=I_q(\widetilde f;0,1)
$$
and
$$
R_q(f;a,b)=(b-a)(1-q)\sum_{k=0}^{\infty}f(a+(b-a)q^k)q^k
=(b-a)I_q(\widehat f;0,1).\ \square
$$
The mentioned connections can be used to derive  the inequalities
for all types of the $q$-integrals. By (\ref{(2.5)}), the
inequalities for the infinite sum $I_q(f;0,b)$ can be derived in
the limit process from this ones for $G_q(f;a,b)$ which are
defined by the finite sum. Using (\ref{(2.6)}) and (\ref{(2.7)}),
the integrals $I_q(f;a,b)$ and $R_q(f;a,b)$ can be considered as
the $q$-integrals over $[0,1]$. Nevertheless, the results for
$I_q(f;a,b)$  are quite rough because the points outside of the
interval of  integration (i.e. points on $[0,a]$) are included.

According to  (\ref{(2.10)}) and Lemma \ref{2.1}, the
following integral relations are valid:
\begin{eqnarray}
&&R_q(fg;a,b)=(b-a)I_q\bigl(\widehat{(fg)};0,1\bigr)
=(b-a)I_q\bigl(\widehat f\ \widehat g;0,1\bigr), \label{(2.8)} \\
&&I_q(fg;a,b)=I_q\bigl(\widetilde{(fg)};0,1\bigr)
=\frac1{b-a}\Bigl(I_q\bigl(\widetilde f \ \widetilde g;0,1\bigr)
-ab\ I_q\bigl(\breve f\ \breve g;0,1\bigr)\Bigr).\qquad \label{(2.9)}
\end{eqnarray}

\section{$q$-Chebyshev inequality}

In this section we give the $q$-analogues of Chebyshev inequality
for the monotonic functions (see \cite{Mitrinovic}, pp. 239.). The
discrete case of this inequality is used in \cite{Gauchman} for
the restricted $q$-integrals. We derive its variants for the rest
of the $q$-integrals.

The function $f:[a,b]\to \mathbb R$ is called $q$-{\it increasing}
($q$-decreasing) on $[a,b]$ if $f(qx)\le f(x)$ ($f(qx)\ge f(x)$)
whenever $x,qx\in[a,b]$. It is easy to see that if the function
$f$ is increasing (decreasing), then it is $q$-increasing
($q$-decreasing) too.

\begin{thm}\label{3.1}
Let $f,g:E_{(J)}\to\mathbb R$ be two real
functions, both $q$-decreasing or both $q$-increasing. If
$J_q(\,\cdot\,;a_{(J)},b)$ is the $q$-integral defined by
$(\ref{(1.1)})$,\ $(\ref{(2.3)})$\ or\ $(\ref{(2.4)})$, it holds
$$
J_q(fg;a_{(J)},b) \ge \frac1{b-a_{(J)}}J_q(f;a_{(J)},b)\
J_q(g;a_{(J)},b). \
$$
\end{thm}

\noindent{\it Proof.} For
$J_q(\,\cdot\,;a_{(J)},b)=G_q(\,\cdot\,;a,b)$,\ $a=bq^n$, the
inequality is proven in \cite{Gauchman}. So, the inequalities
$$
G_q(fg;bq^n,b)\ge\frac 1{b-bq^n}G_q(f;bq^n,b) \ G_q(g;bq^n,b)
$$
are valid for all \ $n=1,2,\ldots$ . When  $n\to\infty$, using
(\ref{(2.5)}) we get the desired inequality for
$J_q(\,\cdot\,;a_{(J)},b)=I_q(\,\cdot\,;0,b)$. In the case
$J_q(\,\cdot\,;a_{(J)},b)=R_q(\,\cdot\,;a,b)$, from the
$q$-monotonicity of the functions $f$ and $g$ on $[a,b]$ follows
the $q$-monotonicity of the functions $\widehat f$ and $\widehat
g$ on $[0,1]$. Hence, we have
$$
I_q(\widehat f \ \widehat g;0,1)\ge I_q(\widehat f;0,1) \
I_q(\widehat g;0,1).
$$
According to (\ref{(2.6)}) and (\ref{(2.7)}) we get the required
inequality. \ $\square$

The Chebyshev inequality in the source form is not valid for
$I_q(\,\cdot\,;a,b)$, where $0<a<b$.

\smallskip

\noindent{\bf Example 3.1} For $f(x)=x^3$ and $g(x)=x^4$ on the
interval $[1,2]$ we have
$$
I_q(x^3\cdot x^4;1,2) - I_q(x^3;1,2)I_q(x^4;1,2)
=255\frac{1-q}{1-q^8} - 465  \frac{(1-q)^2}{(1-q^4)(1-q^5)},
$$
wherefrom we conclude that the inequality holds only for $q>1/2$,
but it has opposite sign for $q<1/2$.

\begin{lem}\label{3.2}
Let the function $f:[0,b]\to\mathbb R$ be increasing and $0<a<b$.
If there exist two positive constants \ $l$ and $L$ such that
$a^2/b^2\le l/L$ and for every $x,y\in [0,b]$ the inequality
$$
l \le \frac{f(x)-f(y)}{x-y} \le L
$$
is valid, then the function $\widetilde f:[0,1]\to\mathbb R$ is
increasing too.
\end{lem}

\noindent{\it Proof.}  Under the conditions of the Lemma, for
every $0\le x <y \le b$ we have
$$
l(y-x)\le f(y)-f(x) \le L(y-x).
$$
Then it holds
$$
\aligned \widetilde f(y)-\widetilde f(x)
&=b\bigl(f(by)-f(bx)\bigr)-a\bigl(f(ay)-f(ax)\bigr)\\
&\ge (b^2 l-a^2 L)(y-x) \ge 0.\ \square
\endaligned
$$

\begin{thm}
Let $f,g:[0,b]\to \mathbb R$ be two real increasing functions. If
there exist the  constants \ $l_f$,\ $L_f$,\ $l_g$\  and\ $L_g$
such that \ $a^2/b^2\le l_f/L_f$,\ \ \ $a^2/b^2\le l_g/L_g$\ \
and
$$
l_f \le \frac{f(x)-f(y)}{x-y}\le L_f,\quad l_g \le
\frac{g(x)-g(y)}{x-y}\le L_g
$$
holds, then the inequalities are valid:
$$
\aligned
(a)\quad  I_q(fg;a,b) & \ge \frac 1{b-a}I_q(f;a,b)
I_q(g;a,b)
-\frac{ab(b-a)}{[3]_q} L_f L_g    \\
(b) \quad I_q(fg;a,b) & \ge \frac 1{b-a}I_q(f;a,b) I_q(g;a,b)
-\frac{ab}{b-a}\bigl(f(b)-f(0)\bigr)\bigl(g(b)-g(0)\bigr).
\endaligned
$$
\end{thm}

\noindent{\it Proof.} Suppose that $f$ and $g$ are both increasing
on $[0,b]$. Then, according to Lemma \ref{3.2}, $\widetilde f$ and
$\widetilde g$ are both increasing and hence $q$-increasing on
$[0,1]$. With respect to (\ref{(2.9)}) we can write
$$
I_q(fg;a,b)=\frac 1{b-a}\Bigl( I_q(\widetilde f\ \widetilde
g;0,1)- ab\ I_q(\breve f\ \breve g;0,1)\Bigr).
$$
Using Theorem \ref{3.1}, we have
$$
I_q(\widetilde f \ \widetilde g;0,1)\ge I_q(\widetilde f;0,1) \
I_q(\widetilde g;0,1),
$$
wherefrom
\begin{equation}
I_q(fg;a,b)\ge \frac 1{b-a}\Bigl(I_q(f;a,b)\
I_q(g;a,b)-abI_q(\breve f\ \breve g;0,1)\Bigr). \label{3.3}
\end{equation}
(a) Under the conditions satisfied by the functions $f$ and $g$ on
$[0,b]$, it holds
$$
\aligned I_q(\breve f\ \breve g;0,1)
&=(1-q)\sum_{k=0}^\infty \bigl(f(bq^k)-f(aq^k)\bigr)\bigl(g(bq^k)-g(aq^k)\bigr)q^k\\
&\le(1-q)\sum_{k=0}^\infty L_f L_g(bq^k-aq^k)^2 q^k
%\\&=(1-q)L_f L_g (b-a)^2 \sum_{k=0}^\infty q^{3k}
=L_f L_g(b-a)^2\frac{1-q}{1-q^3}
\endaligned
$$
Substituting this estimation in (\ref{3.3}), we get the first
inequality.

\noindent(b) Since the functions $f$ and $g$ are increasing on
$[0,b]$, it holds
$$
I_q(\breve f\ \breve g;0,1)\le
(1-q)\bigl(f(b)-f(0)\bigr)\bigl(g(b)-g(0)\bigr)\sum_{k=0}^\infty
q^k =\bigl(f(b)-f(0)\bigr)\bigl(g(b)-g(0)\bigr),
$$
what with (\ref{3.3}) gives the second inequality.\ $\square$

\section{$q$-Gr\"uss inequality}

The  Gr\"uss inequality (see \cite{Mitrinovic}, pp. 296) can be
understood as conversion of Chebyshev one.

\begin{thm}\label{4.1}
Let $f,g:E_{(J)}\to\mathbb R$ be two real
functions, such that $m \le f(x) \le M$,\ $\varphi \le g(x) \le \Phi$ on $E_{(J)}$,
where $m,M,\varphi,\Phi$ are given real constants. If
$J_q(\,\cdot\,;a_{(J)},b)$ is the $q$-integral defined by
$(\ref{(1.1)})$,\ $(\ref{(2.3)})$\ or\ $(\ref{(2.4)})$, it holds
\begin{eqnarray*}
\Biggm| \frac 1{b-a_{(J)}}J_q(fg;a_{(J)},b)-
\frac1{\bigl(b-a_{(J)}\bigr)^2}J_q(f;a_{(J)},b)\ J_q(g;a_{(J)},b)\Biggm| \qquad \qquad \\
\qquad \qquad \le \frac 14(M-m)(\Phi-\varphi).
\end{eqnarray*}
\end{thm}
\noindent{\it Proof.} For the restricted $q$-integrals
$G_q(\,\cdot\,;bq^n,b)$, the inequality is proven in
\cite{Gauchman}. So, for any arbitrary positive integer $n$, the
inequality
\begin{eqnarray*}
\Bigm| \dfrac 1{b-bq^n}G_q(fg;bq^n,b)-\frac
1{(b-bq^n)^2}G_q(f;bq^n,b) \ G_q(g;bq^n,b)\Bigm| \qquad \qquad  \\
 \hspace{3cm}
 \le \frac14(M-m)(\Phi-\varphi)
\end{eqnarray*}
is valid. When $n\to\infty$, we get the required inequality for
$I_q(\,\cdot\,;0,b)$ via (\ref{(2.5)}). \ Finally,  providing the
conditions of the  theorem, the functions $\widehat f$ and
$\widehat g$ are bounded on $[0,1]$ by the constants $m,M,\varphi,
\Phi$ respective. Then,
$$
\Bigm| I_q(\widehat f \ \widehat g;0,1)-I_q(\widehat f;0,1) \
I_q(\widehat g;0,1)\Bigm| \le \frac 14(M-m)(\Phi-\varphi)
$$
holds and  using the relation (\ref{(2.7)}),  we get the
inequality for $R_q(\,\cdot\,;a,b)$.\ \ $\square$

\smallskip

\noindent{\bf Example 4.1} For $f(x)=x$ and $g(x)=x^2$ on the
interval $[1,2]$ we have
$$
I_q(x\cdot x^2;1,2) - I_q(x;1,2)I_q(x^2;1,2)
=(1-2q)\frac{3\ (2-q)}{(1+q)(1+q^2)(1+q+q^2)}.
$$
Including the boundaries of the functions $f(x)$ and $g(x)$, we
can see that the formula of Gr\"uss inequality will not be hold on
for $q\in (0, 1/3)$.

% Theorem 4.4
\begin{thm}
Let $f,g:[0,b]\to \mathbb R$ be two  bounded such that $m\le f(x)
\le M$,\ $\varphi \le g(x) \le \Phi$ on  $[0,b]$, where
$m,M,\varphi,\Phi$ are given real constants. Then it holds
\begin{eqnarray*}
 \Bigm| \frac 1{b-a}I_q(fg;a,b) -
\frac1{(b-a)^2}I_q(f;a,b) I_q(g;a,b)\Bigm| \hspace{4cm}\\
\hspace{5cm}\le \frac
14(M-m)(\Phi-\varphi)\biggl(1+\frac{4ab}{(b-a)^2}\biggr).
\end{eqnarray*}
\end{thm}

\noindent{\it Proof.}  Having in mind the boundaries of $f$ and
$g$ on $[0,b]$, we have
$$
bm-aM\le \widetilde f(x)\le bM-am,\qquad
b\varphi-a\Phi\le\widetilde g(x)\le b\Phi-a\varphi,
$$
where $\widetilde f$ and $\widetilde g$ are the function defined
on $[0,1]$. According to Theorem \ref{4.1}, we have
\begin{eqnarray*}
\Bigm| I_q(\widetilde f\ \widetilde g;0,1)
-I_q(\widetilde f;0,1) \ I_q(\widetilde g;0,1)\Bigm| \hspace{5cm}\\
\hspace{4cm}\le \frac
14(bM-am-bm+aM)(b\Phi-a\varphi-b\varphi+a\Phi).
\end{eqnarray*}
By using (\ref{(2.9)}), we obtain
\begin{eqnarray*}
\phantom{\le} &\Bigm|(b-a)I_q(fg;a,b)-I_q(f;a,b)  \
I_q(g;a,b)\Bigm| -ab\Bigm|I_q(\breve f\ \breve g;0,1)\Bigm|
\hspace{2truecm}
\\
&\hspace{2truecm}\le \Bigm|(b-a)I_q(fg;a,b)-I_q(f;a,b) \
I_q(g;a,b)+ab\ I_q(\breve f\ \breve g;0,1)\Bigm|
\\
 &\hspace{7truecm}\le \frac
14(b-a)^2(M-m)(\Phi-\varphi).
\end{eqnarray*}

With respect to the boundaries of $f$ and $g$ on $[0,b]$, the
estimation
$$
\Bigm|I_q(\breve f\ \breve g;0,1)\Bigm|\le (M-m)(\Phi-\varphi)
$$
holds, what, finally, proves the statement.\ \ $\square$

% Section 5
\section{$q$-Hermite--Hadamard inequality}

The Hermite--Hadamard inequality (see \cite{Mitrinovic}, pp. 10)
is related to the Jensen inequality for the convex function. In
\cite{Gauchman} there is proved a variant of its analogue for the
restricted $q$-integrals. Here we will formulate and prove another
variant of the $q$-Hermite--Hada\-mard inequality for the
restricted $q$-integrals and for the other types of $q$-integrals.

\begin{thm}\label{5.1}
Let $f:[a,b]\to \mathbb{R}$ $(a=bq^n)$ be a convex function. Then
it holds
$$
f\biggl(\frac{a+b}{[2]_q}\biggr) \le \frac1{b-a}G_q(f;a,b) \le
\frac1{[2]_q}\biggl(q\ f\Bigl(\frac aq\Bigr)+f(b)\biggr).
$$
\end{thm}

\noindent{\it Proof.} According to the definition of the
restricted $q$-integral, we have
$$
\frac1{b-a}G_q(f;a,b)=\frac{1-q}{1-q^n} \sum_{k=0}^{n-1}f(bq^k)q^k
=\Biggl(\sum\limits_{k=0}^{n-1}q^k\Biggr)^{-1}
\Biggl(\sum\limits_{k=0}^{n-1}f(bq^k)q^k\Biggr)\ .
$$
If we assign
$$
\overline x=\Biggl(\sum\limits_{k=0}^{n-1}\ q^k\Biggr)^{-1}
\Biggl(\sum\limits_{k=0}^{n-1}bq^k\ q^k\Biggr)
=\frac{b(1+q^n)}{1+q}=\frac{a+b}{1+q}
$$
and apply Jensen inequality for the convex functions
%\cite{Mitrinovic}
on the last term, we obtain
$$
\frac1{b-a}G_q(f;a,b)\ge f(\overline
x)=f\biggl(\frac{a+b}{1+q}\biggr).
$$
On the other side, using a variant of the reverse Jensen
inequality (see \cite{Mitrinovic}, pp. 9.), we get
$$
\aligned \frac1{b-a}G_q(f;a,b)&\le \frac{b-\overline
x}{b-bq^{n-1}}\ f(bq^{n-1})
+\frac{\overline x-bq^{n-1}}{b-bq^{n-1}}\ f(b)\\
&=\biggl(b-\frac aq\biggr)^{-1}\biggl(\Bigl(b-\frac{a+b}{1+q}\Bigr)f\Bigl(\frac
aq\Bigr)+
\Bigl(\frac{a+b}{1+q}-\frac aq\Bigr)f(b)\biggr)\\
&=\frac 1{1+q}\biggl(q\ f\Bigl(\frac aq\Bigr)+f(b)\biggr).\
\square
\endaligned
$$

\begin{thm}\label{5.2}
Let $f:[0,b]\to \mathbb{R}$  be a continuous convex function.
Then,
$$
f\biggl(\frac{b}{[2]_q}\biggr) \le \frac1{b}I_q(f;0,b) \le
\frac1{[2]_q}\Bigl(q\ f(0)+f(b)\Bigr).
$$
\end{thm}

\noindent{\it Proof.}  Since the function $f$ satisfies the
conditions of Theorem \ref{5.1}  on the intervals $[bq^n,b]$ for
every $n\in\mathbb N$, the inequalities
$$
f\biggl(\frac{bq^n+b}{[2]_q}\biggr) \le
\frac1{b-bq^n}G_q(f;bq^n,b) \le \frac1{[2]_q}\biggl(q\
f\biggl(\frac {bq^n}q\biggr)+f(b)\biggr)
$$
are valid. When $n\to\infty$, we obtain the desired inequality
because $f$ is continuous and (\ref{(2.5)}) is satisfied. \
$\square$

\begin{thm}\label{5.3}
Let $f:[a,b]\to \mathbb{R}$ be a continuous convex
function. Then,
$$
f\biggl(\frac{aq+b}{[2]_q}\biggr) \le \frac1{b-a}R_q(f;a,b)) \le
\frac1{[2]_q}\Bigl(q\ f(a)+f(b)\Bigr).
$$
\end{thm}

\noindent{\it Proof.} Under the conditions which are satisfied by
the function $f$ on $[a,b]$, the function $\widehat
f(x)=f\bigl(a+(b-a)x\bigr)$ satisfies the conditions of the
Theorem \ref{5.2} on $[0,1]$. Hence
$$
\widehat f\biggl(\frac{1}{[2]_q}\biggr) \le I_q(\widehat f;0,1)
\le \frac1{[2]_q}\Bigl(q\ \widehat f(0)+\widehat f(1)\Bigr).
$$
According to (\ref{(2.8)}) and the continuity of the function $f$,
we get the desired inequality.\ \ $\square$

Let us remember that the function $f$ is convex on $[0,b]$ if for
all $x,y\in[0,b]$ and $p_1+p_2>0$
$$
f\biggl(\frac{p_1 x+p_2 y}{p_1+p_2}\biggr) \le \frac{p_1 f(x)+p_2 f(y)}{p_1+p_2}
$$
holds. The convexity of the function $\widetilde f$ on $[0,1]$ is due to the existence
of the appropriate constants $l$ and $L$ \ such that the condition
\begin{equation}\label{5.10}
l \le \frac{p_1 f(x)+p_2 f(y)}{p_1+p_2}
-f\biggl(\frac{p_1 x+p_2 y}{p_1+p_2}\biggr) \le L
\end{equation}
is satisfied.

\begin{lem}\label{5.4}
Let the function $f:[0,b]\to\mathbb R$ be convex. If there exist
two positive constants \ $l$ and $L$ such that \ \ $bl \ge aL$ and
for every $x,y\in [0,b]$ and $p_1+p_2>0$\ the condition
$(\ref{5.10})$ is satisfied, then the function $\widetilde
f:[0,1]\to\mathbb R$
%defined by \ $\widetilde f(x)=bf(bx)-af(ax)$
is convex too.
\end{lem}

\noindent{\it Proof.}  Under the conditions of the Lemma, for
every $0\le x,y \le b$ and $p_1+p_2>0$ we have
\begin{eqnarray*}
&&\frac{p_1 \widetilde f(x)+p_2 \widetilde f(y)}{p_1+p_2}  -
\widetilde f\biggl(\frac{p_1 x+p_2 y}{p_1+p_2}\biggr) \\
&&\hspace{4cm}=b\Biggl(\frac{p_1 f(bx)+p_2 f(by)}{p_1+p_2}
-f\biggl(\frac{p_1 bx+p_2 by}{p_1+p_2}\biggr)
\Biggr)\\
&&\hspace{4cm}-a\Biggl(\frac{p_1 f(ax)+p_2 f(ay)}{p_1+p_2}
-f\biggl(\frac{p_1 ax+p_2 ay}{p_1+p_2}\biggr)
\Biggr)\\
&&\hspace{4cm}\ge bl-aL \ge 0.\ \square
\end{eqnarray*}

\begin{thm}\label{5.5}
Let $f:[0,b]\to \mathbb{R}$ be a continuous and convex function.
If there exist two positive constants \ $l$ and $L$ such that \ \
$bl \ge a L$ \ and for every $x,y\in [0,b]$,\ $p_1+p_2>0$ the
condition $(\ref{5.10})$ is satisfied, then it holds
\begin{equation}\label{5.6}
b f\Bigl(\dfrac{b}{[2]_q}\Bigr)-a f\Bigl(\dfrac{a}{[2]_q}\Bigr)\le
I_q(f;a,b) \le \frac{(b-a)qf(0)+ b f(b)-af(a)}{[2]_q}.
\end{equation}
\end{thm}

\noindent{\it Proof.}  According to Lemma \ref{5.4}, the function
$\widetilde f$ is convex on $[0,1]$. Then, using Theorem
\ref{5.2}, we have
$$
\widetilde f\biggl(\frac{1}{[2]_q}\biggr) \le I_q(\widetilde
f;0,1) \le \frac1{[2]_q}\biggl(q\ \widetilde f(0)+\widetilde
f(1)\biggr).
$$
Applying the relation (\ref{(2.6)}) we get the statement. \
$\square$

\begin{cor}\label{5.5}
Let $f:[0,a+b]\to \mathbb{R}$ be a continuous and convex function.
If there exist two positive constants \ $l$ and $L$ such that \ \
$bl \ge a L$ \ and for every $x,y\in [0,a+b]$,\ $p_1+p_2>0$ the
condition $(\ref{5.10})$ is satisfied, then it holds
$$
l+f\Bigl(\frac{a+b}{[2]_q}\Bigr)\le\frac1{b-a} I_q(f;a,b) \le
\frac1{[2]_q}\Bigl(qf(0)+ f(a+b)+L\Bigr).
$$
\end{cor}

\noindent{\it Proof.}  Let $p_1=b/(b-a)$,\ $p_2=-a/(b-a)$.
Applying the condition (\ref{5.10}) with $x=b/(1+q)$,\ $y=a/(1+q)$
on the left term and $x=a$,\ $y=b$ on the right term in
(\ref{5.6}), we get the statement.\ $\square$

% Section 6
\section{The other inequalities}

In this section we will formulate some new inequalities for
$G_q(\,\cdot\,;a,b)$, $I_q(\,\cdot\,;0,b)$ and
$R_q(\,\cdot\,;a,b)$. They will be proven only for
$G_q(\,\cdot\,;a,b)$. In the way presented in the previous
sections, these inequalities for the other two types follow
directly. Furthermore, it seems that the corresponding
inequalities for the integral $I_q(\,\cdot\,;a,b)$ defined by
(\ref{(1.2)}), exist and have different forms because of the
previously mentioned difficulties related to estimating of the
difference of series.

 So, let
$J_q(\,\cdot\,)=J_q(\,\cdot\,;a_{(J)},b)$ denotes the
$q$--integral defined by (\ref{(1.1)}),\ (\ref{(2.3)}) or
(\ref{(2.4)}). In the formulation and proofs of the theorems we
follow the inequalities for the finite sums given in \cite{S.S.
Dragomir}.

The first class are the inequalities the
Cauchy-Buniakowsky-Schwarz type.

% Theorem 6.1
\begin{thm}
Let $f,g:E_{(J)}\to \mathbb R$ be two real functions and $\alpha,
\beta >1$  the numbers satisfying $\dfrac 1\alpha+\dfrac
1\beta=1$. Then the following inequalities hold:

\begin{eqnarray*}
(i) \quad &&\frac 1\alpha J_q(|f|^\alpha) + \frac 1\beta
J_q(|g|^\beta)  \ge \frac 1{b-a_{(J)}} J_q(|f|) J_q(|g|), \\
\\
(ii) \quad &&\frac 1\alpha J_q(|f|^\alpha) J_q(|g|^\alpha) + \frac
1\beta J_q(|f|^\beta) J_q(|g|^\beta)
\ge \Bigl(J_q(|fg|)\Bigr)^2 ,\\
\\
(iii) \quad &&\frac 1\alpha J_q(|f|^\alpha) J_q(|g|^\beta) +\frac
1\beta J_q(|f|^\beta) J_q(|g|^\alpha)
\ge J_q(|f||g|^{\alpha-1})J_q(|f||g|^{\beta-1}),\\
\\
(iv) \quad &&J_q(|f|^\alpha) J_q(|g|^\beta) \ge J_q(|fg|)
J(|f|^{\alpha-1}|g|^{\beta-1}).
\end{eqnarray*}
\end{thm}

\noindent{\it Proof.} If in well-known Young inequality (see
\cite{Mitrinovic}, pp. 381)
$$
\frac 1\alpha x^\alpha+\frac 1\beta y^\beta\ge xy
\qquad\qquad
(x,y\ge0,\quad \alpha, \beta >1:\ \frac 1\alpha+\frac
1\beta=1),
$$
we put $x=|f(bq^i)|$,\ $y=|g(bq^j)|$,\ where $i,j=0,1,\ldots, n-1$, we
have
$$
\frac 1\alpha |f(bq^i)|^\alpha+\frac 1\beta |g(bq^j)|^\beta \ge
|f(bq^i)||g(bq^j)|,\quad i,j=0,1,\ldots, n-1.
$$
Multiplying by $q^{i+j}$ and summing over
$i$ and $j$, we obtain
\begin{eqnarray*}
\frac 1\alpha \sum_{j=0}^{n-1}q^j
\sum_{i=0}^{n-1}q^i|f(bq^i)|^\alpha + \frac 1\beta
\sum_{i=0}^{n-1}q^i \sum_{j=0}^{n-1}q^j|g(bq^j)|^\beta
\ge
\sum_{i=0}^{n-1}q^i|f(bq^i)| \sum_{j=0}^{n-1}q^j|g(bq^j)|
\end{eqnarray*}
and, finally, inequality ({\it i}). The rest of inequalities can
be proved in the same manner by the next choice of the parameters
in Young inequality:
\begin{eqnarray*}
(ii)\ & x=|f(bq^j)|\ |g(bq^i)|,  \  &y=|f(bq^i)|\ |g(bq^j)|,\\
(iii)\  & x=|f(bq^j)|/|g(bq^j)|, \ &y=|f(bq^i)|/|g(bq^i)|,\quad \bigl(g(bq^j) \ g(bq^i)\ne 0\bigr),\\
(iv)\  & x=|f(bq^i)|/|f(bq^j)|,  \ &y=|g(bq^i)|/|g(bq^j)|,\quad
\bigl(f(bq^j)\ g(bq^j)\ne 0\bigr),
\end{eqnarray*}
where additional conditions about not vanishing for $f$ and $g$ do
not have influence on final conclusion. $\square$

% Theorem 6.2
\begin{thm}
Let $f,g:E_{(J)}\to \mathbb R$ be two real functions and $\alpha,
\beta >1$  the numbers satisfying $\dfrac 1\alpha+\dfrac
1\beta=1$. Then the following inequalities hold:
$$
\aligned
(i) \ &\frac 1\alpha J_q(|f|^\alpha) J_q(|g|^2) +
\frac 1\beta J_q(|f|^2) J_q(|g|^\beta)
\ge J_q(|fg|) J_q(|f|^{2/\beta}|g|^{2/\alpha}),\\
(ii) \ &\frac 1\alpha J_q(|f|^2) J_q(|g|^\beta) +
\frac 1\beta J_q(|f|^\alpha) J_q(|g|^2)
\ge J_q(|f|^{2/\alpha}|g|^{2/\beta})
J_q(|f|^{\alpha-1}|g|^{\beta-1}), \\
(iii) \ & \qquad\qquad \ J_q(|f|^2) J_q\Bigl(\frac 1\alpha |g|^\alpha
+\frac 1\beta |g|^\beta\Bigr)\ge J_q(|f|^{2/\alpha}|g|)
J_q(|f|^{2/\beta}|g|).
\endaligned
$$
\end{thm}

\noindent{\it Proof.}
As previous, the proof is based on Young inequality with appropriate
choice of the parameters:

\begin{eqnarray*}
(i)\ &x=|f(bq^i)|\ |g(bq^j)|^{2/\alpha}, \
&y=|f(bq^j)|^{2/\beta}\ |g(bq^i)|,\\
(ii)\ &x=|f(bq^i)|^{2/\alpha}/|f(bq^j)|, \
&y=|g(bq^i)|^{2/\beta}/|g(bq^j)|\quad
\ (f(bq^j)g(bq^j)\ne 0),\\
(iii)\ &x=|f(bq^i)|^{2/\alpha}\ |g(bq^j)|, \
&y=|f(bq^j)|^{2/\beta}\ |g(bq^i)|.\ \square
\end{eqnarray*}

The following few inequalities include  the boundaries of the functions.

% Theorem 6.3
\begin{thm}\label{6.3}
If $f,g:E_{(J)}\to \mathbb R$ are two positive functions and
$$
m=\min_{a\le x\le b}\frac{f(x)}{g(x)},\qquad M=\max_{a\le x\le
b}\frac{f(x)}{g(x)},
$$
then the following inequalities hold:
$$
\aligned
(i) \quad & \qquad \qquad \quad 0\le J_q(f^2)J_q(g^2) \le
\frac{(m+M)^2}{4mM}\Bigl(J_q(fg)\Bigr)^2,\\
(ii)\quad & 0\le \sqrt{(J_q(f^2) J_q(g^2)}-J_q(fg)
\le \frac{(\sqrt{M}-\sqrt{m})^2}{2\sqrt{mM}}J_q(fg),\\
(iii)\quad & 0\le J_q(f^2) J_q(g^2)-\Bigl(J_q(fg)\Bigr)^2
\le \frac{(M-m)^2}{4mM}\Bigl(J_q(fg)\Bigr)^2.
\endaligned
$$
\end{thm}

\noindent{\it Proof.} With respect to the definition of
$G_q(\,\cdot\,;a,b)$, the inequality ({\it i}) is the immediate
consequence of the Cassels inequality (see \cite{S.S. Dragomir},
pp. 72). The inequalities ({\it ii}) and ({\it iii}) can be
obtained by a few transformations of ({\it i}).\ $\square$

% Theorem 6.4
\begin{thm}
If $f,g:E_{(J)}\to \mathbb R$ are two positive functions such that
$$
0<c\le f(x)\le C<\infty,\qquad 0<d\le g(x)\le D<\infty,
$$
then the following inequalities hold:
$$
\aligned
(i)\quad \ & \qquad  \qquad  \qquad  0\le J_q(f^2)J_q(g^2) \le
\frac{(cd+CD)^2}{4cdCD}\Bigl(J_q(fg)\Bigr)^2,\\
(ii)\quad  & \quad 0\le \sqrt{J_q(f^2) J_q(g^2)}-J_q(fg)
\le \frac{(\sqrt{CD}-\sqrt{cd})^2}{2\sqrt{cdCD}}J_q(fg),\\
(iii)\quad  & 0\le J_q(f^2) J_q(g^2)-\Bigl(J_q(fg)\Bigr)^2
\le \frac{(CD-cd)^2}{4cdCD}\Bigl(J_q(fg)\Bigr)^2.
\endaligned
$$
\end{thm}

\noindent{\it Proof.}
Under the conditions satisfied by the functions $f$ and $g$, we have
$$
\frac  cD \le \frac{f(x)}{g(x)} \le \frac Cd
%,\quad a\le x\le b
.
$$
Applying Theorem \ref{6.3}  we get the inequality ({\it i}) and, using
it, ({\it ii}) and ({\it iii}).\ $\square$

% Theorem 6.5
\begin{cor}
Let $f:E_{(J)}\to \mathbb R$ be a positive function such that
$$
0<c\le f(x)\le C<\infty.
$$
Then the following inequality holds:
$$
J_q(f^2) \le \frac{(c+C)^2}{4cC\bigl(b-a_{(J)}\bigr)}\Bigl(J_q(f)\Bigr)^2.
$$
\end{cor}

The next few inequalities are obtained via Jensen inequality for
the convex functions.

% Theorem 6.6
\begin{thm}
Let $f,g:E_{(J)}\to \mathbb R$ be two positive functions and $p\ne
0$ a real number. Then it holds
$$
\aligned
\Bigl(J_q(fg)\Bigr)^p &\le \Bigl(J_q(f^2)\Bigr)^{p-1}
J_q(f^{2-p}g^p),\quad \textrm{for} \quad p\notin (0,1),\\
\Bigl(J_q(fg)\Bigr)^p &\ge \Bigl(J_q(f^2)\Bigr)^{p-1}
J_q(f^{2-p}g^p),\quad \textrm{for} \quad p\in (0,1).
\endaligned
$$
\end{thm}

\noindent{\it Proof.}  For $p\notin (0,1)$ the function $t\mapsto
t^p$ is convex. Applying the Jensen inequality for convex
functions (see \cite{Mitrinovic}, pp.6.) we have
$$
\Biggl(\frac{\sum\limits_{k=0}^{n-1} f(bq^k) g(bq^k)q^k}
{\sum\limits_{k=0}^{n-1}\bigl(f(bq^k)\bigr)^2 q^k}\Biggr)^p \le
\frac 1{\sum\limits_{k=0}^{n-1}\bigl(f(bq^k)\bigr)^2 q^k} \
\sum_{k=0}^{n-1} \Bigl(\frac{g(bq^k)}{f(bq^k)}\Bigr)^p
\bigl(f(bq^k)\bigr)^2 q^k,
$$
i.e.,
$$
\biggl(\sum_{k=0}^{n-1} f(bq^k) g(bq^k)q^k\biggr)^p \le
\biggl(\sum_{k=0}^{n-1}\bigl(f(bq^k)\bigr)^2 q^k\biggr)^{p-1}
\biggl(\sum_{k=0}^{n-1} \bigl(g(bq^k)\bigr)^p
\bigl(f(bq^k)\bigr)^{2-p} q^k\biggr).
$$
According to the definition of $G_q(\,\cdot \,;a,b)$
we get the inequality. The reverse case is obtained
for $p\in (0,1)$ because of the concave function $t\mapsto t^p$. \
$\square$

% Theorem 6.7
\begin{cor}
Let $f:E_{(J)}\to \mathbb R$ be a  positive function and $p\ne 0$
a real number. Then it holds
$$
\Bigl(J_q(f)\Bigr)^p \le \bigl(b-a_{(J)}\bigr)^{p-1}\  J_q(f^p),
$$
for $p\notin (0,1)$, or reverse for $p\in(0,1)$.
\end{cor}

% Theorem 6.8
\begin{thm}
If $f,g:E_{(J)}\to \mathbb R$ are two positive functions such that
$$
0<m\le \frac{g(x)}{f(x)}\le M<\infty
$$
and $p\ne 0$ a real number, then it holds
$$
\aligned
J_q(f^{2-p}g^p)+\frac{mM(M^{p-1}-m^{p-1})}{M-m}
J_q(f^p) \le \frac{M^p-m^p}{M-m}J_q(fg),
\endaligned
$$
for $p\notin (0,1)$, or reverse for $p\in(0,1)$. Especially, for
$p=2$, we have
$$
J_q(g^2) + mM J_q(f^2)\le (M+m) J_q(fg).
$$
\end{thm}

\noindent{\it Proof.} The inequality is based on the Lah-Ribari\'c
inequality (see \cite{Mitrinovic}, pp. 9., \cite{S.S. Dragomir},
pp. 123).\ $\square$

\medskip

\centerline{\bf Acknowledgement }
\smallskip

This research was supported by the Science Foundation of Republic
Serbia, Project No. 144023 and Project No. 144013.

\end{document}